\documentclass[12pt,a4paper,twoside]{article}

\usepackage{amsmath,amsfonts,amsthm,amsopn,color,amssymb,cite}
\usepackage{palatino}
\usepackage{graphicx}

\newcommand{\R}{\mathbb{R}}

\newcommand{\RN}{{\mathbb{R}^N}}
\newcommand{\RT}{{\mathbb{R}^3}}

\renewcommand{\le}{\leqslant}
\renewcommand{\ge}{\geqslant}
\renewcommand{\a }{\alpha }

\renewcommand{\d }{\delta }

\renewcommand{\l }{\lambda}
\newcommand{\n }{\nabla }

\renewcommand{\t}{\theta}
\renewcommand{\O}{\Omega}

\newcommand{\G}{\Gamma}

\newcommand{\Ne}{\mathcal{N}}
\newcommand{\E}{\mathcal{E}}

\newcommand{\N}{\mathbb{N}}

\newcommand{\D }{{\mathcal D}^{1,2}(\RT)}
\renewcommand{\H}{H^{1}(\RT)}
\newcommand{\irn }{\int_{\RN}}
\newcommand{\irt }{\int_{\RT}}

\def\bbm[#1]{\mbox{\boldmath $#1$}}

\newtheorem{theorem}{Theorem}[section]
\newtheorem{lemma}[theorem]{Lemma}

\newtheorem{remark}[theorem]{Remark}

\renewenvironment{proof}{\noindent{\textbf{Proof\quad}}}{$\hfill\square$\vspace{0.2 cm}\\}
\newenvironment{proofmain}{\noindent{\textbf{Proof of Theorem  \ref{main}\quad}}}{$\hfill\square$\vspace{0.2 cm}\\}

\title{{\bf Ground state solutions for the nonlinear Schr\"odinger-Maxwell
equations with a singular potential}}

\author{A. Azzollini \thanks{Dipartimento di Matematica, Universit\`a degli
Studi di Bari,  Via E. Orabona 4, I-70125 Bari, Italy, e-mail:
{\tt azzollini@dm.uniba.it}}
 \; \& \;
A. Pomponio\thanks{Dipartimento di Matematica, Politecnico di
Bari, Via E. Orabona 4, I-70125 Bari, Italy, e-mail: {\tt
a.pomponio@poliba.it}}}

\date{}

\begin{document}

\maketitle
\begin{abstract}
In this paper we find a ground state solution for the nonlinear Schr\"odinger-Maxwell equations
\[
\left\{
\begin{array}{ll}
-\Delta u + V(x)u + \phi u = |u|^{p-1}u  & \hbox{in }\RT,
\\
-\Delta \phi = u^2 & \hbox{in }\RT.
\end{array}
\right.
\]
where $V$ is a possibly singular potential and $3<p<5$.
\end{abstract}

\section{Introduction}
In this paper we study the following problem:
\begin{equation}\label{ScMa}\tag{$\mathcal{SM}$}
\left\{
\begin{array}{ll}
-\Delta u + V(x)u + \phi u = f(u)  & \hbox{in }\RT,
\\
-\Delta \phi = u^2 & \hbox{in }\RT,
\end{array}
\right.
\end{equation}
where $V:\RT \to \R$ and $f\in C(\R,\R).$ Such a system represents
the nonlinear Schr\"odinger-Maxwell equations in the electrostatic
case (for more details on the physical aspects of this problem and
the relativistic Klein-Gordon-Maxwell equations, we refer to
\cite{BF,BF2}). In \cite{BF}, the potential $V$ has been supposed
constant, and the linear version of the problem (i.e. $f \equiv
0$) has been studied as an eigenvalue problem for a bounded
domain. The linear and nonlinear Schr\"odinger-Maxwell equations
have been treated also in \cite{AR, CS,C1,C2,CG,DM,DW1,DW2,DA,K,Ruiz,Ru,S},
where $V$ is a positive constant or a radially symmetric
potential. For a related problem see \cite{PS,RS}.

Recently, the case of a positive and bounded non-radial potential
$V$ has been studied in \cite{WSZ}, when $f$ is asymptotically
linear, and in \cite{AP}, when $f(u)=|u|^{p-1}u$, with $3<p<5$.
Moreover, in \cite{AP}, existence of ground state solutions for
problem \eqref{ScMa} has been proved in several situations,
including the positive constant potential case.

In this paper we are interested in looking for solutions to the
problem
\begin{equation}\label{SM}
\left\{
\begin{array}{ll}
-\Delta u + V(x)u + \phi u = |u|^{p-1}u  & \hbox{in }\RT,
\\
-\Delta \phi = u^2 & \hbox{in }\RT,
\end{array}
\right.
\end{equation}
where $3<p<5$ and $V$ satisfies the following hypotheses:
\begin{itemize}
\item[({\bf V1})] $V\colon\RT \to\R$ is a measurable function;
\item[({\bf V2})] $V_\infty:=\liminf_{|y|\to\infty}V(y) \ge V(x)$, for almost every $x\in \RT$, and the inequality is strict in a non-zero measure domain;
\item[({\bf V3})] there exists $\bar C>0$ such that, for any  $u\in\H,$
\begin{align*}
\irt |\n u|^2 + V(x) u^2 \ge \bar C \|u\|^2.
\end{align*}
\end{itemize}

In the same spirit of \cite{L}, where similar hypotheses on $V$
are introduced to study singular nonlinear Schr\"odinger
equations, our aim is to extend the existence result contained in
\cite{AP} to the case of a potential unbounded from below.

\begin{remark}\label{re:pot}
Hypotheses on $V$ are satisfied by a large class of potentials including those most meaningful by a physical point of view. Here we give some examples of admissible potentials $V\colon \RT \to \R$:
\begin{enumerate}
\item $V(x)=V_1- \l |x|^{-\a}$,  where $V_1$ is a positive constant, $\a=1,2$ and $\l$ is a positive
constant small enough;
\item $V(x)=V_1(x) -\l |x|^{-\a}$,  where
$V_1$ is a potential bounded below by a positive constant and
satisfying ({\bf V2}), $\a=1,2$ and $\l$ is a sufficiently small
positive constant;
\item $V(x)=V_1(x)- \l V_2(x)$,  where $V_1$ is
a potential bounded below by a positive constant and satisfying
({\bf V2}), $\l$ is a sufficiently small positive constant and
$V_2$ is a positive function such that
\[
\exists \a_1 >0 , \a_2\ge 0 \colon \!\!\irt V_2(x) u^2 \le
\!\irt \!\a_1 |\n u|^2 + \a_2 u^2,\; \hbox{ for any } u \in
H^1(\RT),
\]
and
\[
\lim_{|x|\to +\infty}V_2(x)=0.
\]
\end{enumerate}
\end{remark}

The solutions $(u,\phi)\in \H \times \D$ of
\eqref{SM} are the critical points of the action functional $\mathcal{E}
\colon \H \times \D \to \R$, defined as
\[
\mathcal{E}(u,\phi):=\frac 12 \irt |\n u|^2 + V(x)u^2
-\frac 14 \irt |\n \phi|^2
+\frac 12 \irt \phi u^2
-\frac 1{p+1} \irt |u|^{p+1}.
\]
We are looking for a ground state solution of \eqref{SM}, that is,
according to the definition given in \cite{CGM}, a solution
$(u_0,\phi_0)$ of \eqref{SM} with the property of having the least
action among all possible solutions of \eqref{SM}, namely
$\E(u_0,\phi_0)\le \E(u,\phi)$, for any solution $(u,\phi)$ of
\eqref{SM}.

Our main result is the following
\begin{theorem}\label{main}
If $V$ satisfies ({\bf V1-3}), then the problem \eqref{SM} has a ground state
solution for any $p\in ]3,5[.$
\end{theorem}

The action functional $\E$ exhibits a strong indefiniteness, namely it is
unbounded both from below and from above on infinite dimensional
subspaces. This indefiniteness can be removed using the reduction
method described in \cite{BFMP}, by which we are led to study a
one variable functional that does not present such a strongly
indefinite nature.

The main difficulty related with the problem of finding the
critical points of the new functional, consists in the lack of
compactness of the Sobolev spaces embeddings in the unbounded
domain $\RT$. In \cite{AP}, this difficulty has been overcome by
means of a concentration-compactness argument on suitable measures
(see the definition of the sequence $\mu_n$ in \cite[Section
2.3.1]{AP}). Actually such measures, related with a minimizing
sequence for the functional restricted to the Nehari manifold,
seem to be not the right ones in our situation. In fact, since the
potential $V$ is permitted to be unbounded below, we have no way
to affirm that the integral
    $$\int_\O|\n u|^2+ V(x) |u|^2$$
is nonnegative for any $u\in\H$ and $\O\subset \RT.$

This technical difficulty has been overcome exploiting a more
suitable version of the concentration-compactness principle, based
on that used, for example, in \cite[Lemma 6.1]{CL}.

Section \ref{se:pre} is devoted to some preliminaries necessary to prove Theorem~\ref{main} in Section \ref{se:proof}.

\vspace{0.5cm}
\begin{center}
{\bf NOTATION}
\end{center}

\begin{itemize}
\item For any $1\le s< +\infty$, $L^s(\RT)$ is the usual Lebesgue space endowed with the norm
\[
\|u\|_s^s:=\irt |u|^s;
\]
\item $\H$ is the usual Sobolev space endowed with the norm
\[
\|u\|^2:=\irt |\n u|^2+ u^2;
\]
\item $\D$ is completion of $C_0^\infty(\RT)$ with respect to the norm
\[
\|u\|_{\D}^2:=\irt |\n u|^2;
\]
\item for any $r>0$ and $A\subset \RT$
\begin{align*}
B_r    &:=\{y\in\RT\mid |y|\le r\},\\
A^c    &:= \RT\setminus A;
\end{align*}
\item $C,\,C',\,C_i$ are positive constants which can change from
line to line;
\item $o_n(1)$ is a quantity which goes to zero as $n \to +\infty$.
\end{itemize}

\section{Preliminary results}\label{se:pre}

We first recall some well-known facts (see, for instance
\cite{BF,C1,C2,CG,DM,Ru}). For every $u\in L^{12/5}(\RT)$, there
exists a unique $\phi_u\in \D$ solution of
\[
-\Delta \phi=u^2,\qquad \hbox{in }\RT.
\]
It can be proved that $(u,\phi)\in H^1(\RT)\times \D$ is a
solution of \eqref{SM} if and only if $u\in\H$ is a critical point
of the functional $I\colon \H\to \R$ defined as
\begin{equation}\label{eq:defI}
I(u)= \frac 12 \irt |\n u|^2 +  V(x) u^2 +\frac 14 \irt \phi_u u^2
-\frac{1}{p+1}\irt |u|^{p+1},
\end{equation}
and $\phi=\phi_u$.\\
The functions $\phi_u$ possess the following properties (see
\cite{DM} and \cite{Ru})
\begin{lemma}\label{le:prop}
For any $u\in\H$, we have:
\begin{itemize}
\item[i)] $\|\phi_u\|_{\D}\le C \|u\|^2,$ where $C$ does not
depend from $u.$ As a consequence there exists $C'>0$ such that
$$
\irt\phi_u u^2\le C'\|u\|_\frac {12}5^4;
$$
\item[ii)] $\phi_u\ge 0;$ \item[iii)] for any $t>0$:
$\phi_{tu}=t^2\phi_u;$
    \item[iv)] for any $\O\subset\RT$
measurable,
$$
\int_{\O}\phi_u u^2=\int_{\O}\irt\frac{u^2(x)u^2(y)}{|x-y|}dx\,dy.
$$
\end{itemize}
\end{lemma}

In order to get our result, we will use a very standard device: we
will look for a minimizer of the functional \eqref{eq:defI}
restricted to the Nehari manifold
\begin{equation*}
\Ne=\left\{u\in\H\setminus\{0\}\mid G(u)=0\right\},
\end{equation*}
where
\[
G(u):= \irt  |\n u|^2 +   V(x) u^2 +  \phi_u u^2 -
|u|^{p+1}.
\]
The following lemma describes some properties of the Nehari
manifold~$\Ne$ (see \cite{AP,R}):
\begin{lemma}\label{le:N}
\begin{enumerate}
\item For any $u\neq 0$ there exists a unique number $\bar t>0$
such that $\bar t u\in \Ne$ and
\[
I(\bar t u)=\max_{ t \ge 0}I( t u);
\]
\item there exists a positive constant $C$, such that for all $u
\in \Ne$, $\|u\|_{p+1}\ge C$; \item $\Ne$ is a $C^1$ manifold.
\end{enumerate}
\end{lemma}

The Nehari manifold $\Ne$ is a natural constrained for the
functional $I,$ therefore we are allowed to look for critical
points of $I$ restricted to $\Ne$.

In view of this, we assume the following definition
\begin{equation*}
c:=\inf_{u\in\Ne} I(u),
\end{equation*}
so that our goal is to find $\bar u\in\Ne$ such that $I(\bar u)=c$, from which we would deduce that $(\bar u,\phi_{\bar u})$ is a
ground state solution of \eqref{SM}.

First we recall some preliminary lemmas which can be obtained by
using the same arguments as in \cite{R} (see also \cite{AP}).

As a consequence of the Lemma \ref{le:N}, we are allowed to define
the map $t:\H\setminus\{0\}\to\R_+$ such that for any $u\in\H,$
$u\neq 0:$
\begin{equation*}
I\left( t(u)u\right)=\max_{t\ge 0} I(t u).
\end{equation*}
\begin{lemma}\label{le:ccc2}
The following equalities hold
\begin{equation*}
c=\inf_{g\in \G} \max_{t\in [0,1]} I(g(t))=\inf_{u\neq 0} \max_{t
\ge 0} I( tu),
\end{equation*}
where
\begin{equation*}
\G=\left\{g\in C\big([0,1],\H\big) \mid g(0)=0,\;I(g(1))\le 0,
\;g(1)\neq 0\right\}.
\end{equation*}
\end{lemma}

Let $\l>0$, we define
\begin{align*}
I_\l(u) &:=\frac 12 \irt |\n u|^2 + \l u^2 +\frac 14 \irt \phi_u
u^2 -\frac{1}{p+1}\irt |u|^{p+1},
\\
c(\l)    &:=\inf_{u\in \Ne_\l}I_\l(u),
\end{align*}
where $\Ne_\l$ is the Nehari manifold of $I_\l$.

\begin{lemma}\label{le:vn}
Let $\l,\:\l',\:\l_n>0$.
\begin{enumerate}
\item If $\l < \l'$, then $c(\l) < c(\l')$;
\item if $\l_n\to \l$, then $c(\l_n)\to c(\l).$
\end{enumerate}
\end{lemma}

As in \cite{R}, we have
\begin{lemma}\label{le:cinfty}
If $V$ satisfies ({\bf V1-3}), we get $c<c(V_\infty)$.
\end{lemma}

\begin{proof}
By \cite[Theorem 1.1]{AP}, there exists $(w,\phi_w)\in\H\times \D$
a ground state solution of the problem
\begin{equation*}
\left\{
\begin{array}{ll}
-\Delta u + V_\infty u + \phi u = |u|^{p-1} u & \hbox{in }\RT,
\\
-\Delta \phi = u^2 & \hbox{in }\RT.
\end{array}
\right.
\end{equation*}
Let $t(w)>0$ be such that $t(w)w\in \Ne$. By ({\bf V2}), we have
\begin{align*}
c(V_\infty) & =I_{V_\infty}(w) \ge I_{V_\infty} \big(t(w) w\big)
\\
&=I\big(t(w) w\big) + \irn \big(V_\infty - V(x)\big) |t(w)w|^2
> c,
\end{align*}
and then we conclude.
\end{proof}

\section{Proof of the main theorem}\label{se:proof}

Let $(u_n)_n \subset \Ne$ such that
\begin{equation}\label{eq:lim}
\lim_n I(u_n)=c.
\end{equation}
We define the functional $J\colon \H\to\R$ as:
\begin{equation*}
J(u)= \left(\frac {1}{2}-\frac 1 {p+1}\right) \irt|\n u|^2 +
 V(x) u^2 + \left(\frac{1}{4}-\frac 1 {p+1}\right)\irt
\phi_{u} u^2.
\end{equation*}
Observe that for any $u\in\Ne,$ we have $I(u)=J(u).$
\\
By ({\bf V3}) and \eqref{eq:lim}, we deduce that $(u_n)_n$ is
bounded in $\H,$ so there exists $\bar u\in\H$ such that, up to a
subsequence,
\begin{align}
&u_n\rightharpoonup \bar u\quad\hbox{weakly in }\H,
\label{eq:weakM}
\\
&u_n\to \bar u\quad\hbox{in }L^s(B), \hbox{with }B\subset \RT, \hbox{bounded, and }1\le s<6.  \nonumber 
\end{align}
To prove Theorem \ref{main}, we need some compactness on the
sequence $(u_n)_n.$\\
We denote by $\nu_n$ the measure
\begin{equation*}
\nu_n(\O)=  \left(\frac {1}{2}-\frac 1 {p+1}\right) \int_\O |\n
u_n|^2 + V(x) u_n^2 + \left(\frac {1}{4}-\frac 1 {p+1}\right)
\int_\O \phi_{u_n} u_n^2.
\end{equation*}
Observe that, since there is no lower boundedness condition on the
potential $V$, the measures $\nu_n$ may be not positive, and then
we are not allowed to use the Lions' concentration arguments \cite{L1,L2} on
them. However, using a variant presented in \cite{CL}, in the
following theorem we are able to show that the
functions $u_k$ concentrate in the $\H-$norms.

\begin{theorem}\label{th:conc}
For any $\d>0$ there exists $\tilde R>0$ such that for any $n\ge
\tilde R$
\begin{equation*}
\int_{|x|>\tilde  R} (|\n u_n|^2+ |u_n|^2) < \d.
\end{equation*}
\end{theorem}
\begin{proof}
By contradiction, suppose that there exist $\d_0>0$ and a
subsequence $(u_k)_k$ such that for any $k\ge 1$
\begin{equation}\label{eq:contr}
\int_{|x| > k} (|\n u_k|^2+ |u_k|^2) \ge \d_0.
\end{equation}
We define
\begin{equation*}
\rho_k(\O)=  \int_\O
|\n u_k|^2 + |u_k|^2 + \int_\O \phi_{u_k} u_k^2.
\end{equation*}
and, for any $r>0$, we set $A_r:=\{x\in\RT\mid r\le |x| \le
r+1\}.$
\\
We claim that
\begin{equation}\label{eq:rho}
\hbox{for any $\mu>0$ and $R>0$, there exists }r> R\hbox{ such that
}\rho_k(A_r)<\mu
\end{equation}
for infinitely many $k.$ If not, then there should exist
$\hat\mu>0$ and $\widehat R \in \N$  such that, for any $m\ge \widehat
R$, there exists $p(m)$ such that, for any $k\ge p(m)$,
\begin{equation*}
\rho_k(A_m)\ge \hat\mu.
\end{equation*}
We are allowed to take $(p(m))_m$ not decreasing, so that for
every $m\ge \widehat R$ we could get $u_k$ such that, using $i$ of
Lemma \ref{le:prop},
$$
C\|u_k\|^2(1+\|u_k\|^2)\ge\|u_k\|^2 + \irt \phi_{u_k}u_k^2 \ge
\rho_k(B_m\setminus B_{\widehat R}) \ge \left(m-\widehat
R\right)\hat\mu
$$
contradicting the boundedness in $\H$ of the sequence $(u_n)_n.$
\\
So, we assume that \eqref{eq:rho} holds. Taking into account Lemma \ref{le:vn} and Lemma
\ref{le:cinfty}, consider $\mu>0$ such that
\begin{equation*}
c < c(V_\infty-\mu) < c(V_\infty).
\end{equation*}
Using ({\bf V2}), there exists $R_\mu\in \N$ such that for almost every
$|x|\ge R_\mu$
\begin{equation}\label{eq:infty}
V( x)\ge V_\infty - \mu>0;
\end{equation}
we take $r> R_\mu$ such that, up to a subsequence,
\begin{equation}\label{eq:mu}
\rho_k(A_r)<\mu,\quad\hbox{ for all }k\ge 1.
\end{equation}
In particular, \eqref{eq:infty} and \eqref{eq:mu} imply
\begin{align}
\int_{A_r}|\n u_k|^2 +   V(x)u_k^2 = O(\mu),&\quad\hbox{ for all
}k\ge 1, \label{eq:mu1}
\\
\int_{A_r}\phi_{u_k} u_k^2 = O(\mu) ,&\quad\hbox{ for all }k\ge 1.
\label{eq:mu2}
\end{align}
Let $\chi \in C^\infty,$ such that $\chi=1$ in $B_r$ and $\chi=0$
in $(B_{r+1})^c,$ $0\le\chi\le 1$ and $|\n \chi|\le 2.$ Set
$v_k=\chi u_k$ and $w_k=(1-\chi)u_k.$
\\
By simple computations, by \eqref{eq:infty} and \eqref{eq:mu1} we
infer
\begin{align}
\int_{A_r}|\n v_k|^2 +   V(x) v^2_k = O(\mu), \qquad \int_{A_r}| v_k|^{p+1} = O(\mu),\nonumber
\\
\int_{A_r}|\n w_k|^2 +   V(x) w^2_k = O(\mu), \qquad \int_{A_r}| w_k|^{p+1} = O(\mu).\nonumber
\end{align}
Hence, we deduce that
\begin{align}
\irt |\n u_k|^2 +   V(x) u^2_k &=\irt|\n v_k|^2 + V(x)v^2_k
\nonumber
\\
&\quad+\irt|\n w_k|^2 + V(x) w^2_k +O(\mu), \label{eq:norm}
\\
\irt |u_k|^{p+1}&=\irt|v_k|^{p+1}+\irt|w_k|^{p+1}
+O(\mu)\label{eq:p+1};
\end{align}
for large $k\ge 1,$ by  \eqref{eq:contr} and \eqref{eq:infty}, we
also deduce that there exists $\d'>0$ such that
\begin{equation}\label{eq:contr2}
\irt |\n w_k|^2 + V(x) |w_k|^2 \ge \d'+ O(\mu).
\end{equation}
Moreover, by point $iv$ of Lemma \ref{le:prop} and \eqref{eq:mu2},
we have
\begin{align}
\irt \phi_{u_k} u_k^2 &= \irt \phi_{v_k} v_k^2 + \irt \phi_{w_k}
w_k^2 +
2\int_{B_{r}}\!\int_{B_{r+1}^c}\!\!\!\!\frac{u_k^2(x)u_k^2(y)}{|x-y|}d
x\,d y + O(\mu) \nonumber
\\
&\ge \irt \phi_{v_k} v_k^2 + \irt \phi_{w_k} w_k^2 + O(\mu).
\label{eq:phi}
\end{align}
Hence, by \eqref{eq:norm} and \eqref{eq:phi}, we get
\begin{align}
J(u_k) \ge J(v_k)+J(w_k)+O(\mu),\nonumber
\end{align}
and then, using \eqref{eq:contr2} and ({\bf V3}), we deduce
\begin{align}
J(u_k) - C\d' & \ge J(v_k) + O(\mu), \label{eq:vk}
\\
J(u_k)      & \ge J(w_k) + O(\mu).  \label{eq:wk}
\end{align}
We recall the definition of the functional $G\colon \H\to \R$
\[
G(u)=\irt |\n u|^2 +  V(x) u^2+ \phi_u u^2- |u|^{p+1}
\]
and that if $u\in \Ne$, then $G(u)=0$. By \eqref{eq:norm},
\eqref{eq:p+1}  and \eqref{eq:phi}, we have
\begin{equation}
0=G(u_k) \ge G(v_k)+G(w_k)+O(\mu)\label{eq:G}.
\end{equation}
We have to distinguish three cases.
\\
\
\\
{\sc Case 1:} up to a subsequence, $G(v_k) \le 0$.
\\
By Lemma \ref{le:N}, for any $k \ge 1$, there exists $\t_k>0$ such
that $\t_k v_k \in \Ne$, and then
\begin{equation}\label{eq:nehari}
\irt |\n v_k|^2 + V(x) v_k^2 +  \t_k^2 \phi_{v_k}v_k^2 =
\irt\t_k^{p-1} |v_k|^{p+1}.
\end{equation}
By \eqref{eq:nehari} we have
\begin{align*}
(\t_k^{p-1}-1)\irt  |\n v_k|^2 + V(x) v_k^2 +
(\t_k^{p-1}-\t_k^2)\irt\phi_{v_k}v_k^2 \le 0,
\end{align*}
and, by ({\bf V3}), we deduce that $\t_k\le 1$. Therefore, for all $k\ge 1,$ by
({\bf V3}) and  \eqref{eq:vk},
\begin{align*}
c \! \le \! I( \t_k v_k)\! = \!J( \t_k v_k) \!\le\! J(v_k)\! \le
\! J(u_k) -C \d' +O(\mu) =  c -C \d' + o_k(1) + O(\mu),
\end{align*}
which is a contradiction.
\\
\
\\
{\sc Case 2:} up to a subsequence, $G(w_k) \le 0$.
\\
Let $(\eta_k)_k$ be such that, for any $k\ge 1,$ $\eta_k
w_k\in\Ne.$ Arguing as in the previous case, we deduce that
$\eta_k\le 1.$ Define $\tilde w_k=\eta_k w_k.$ Let $(t_k)_k$ be
such that, for any $k\ge 1,$ $t_k\tilde w_k\in \Ne_{V_\infty -\mu}.$
\\
By \eqref{eq:infty},
\begin{multline*}
\irt |\n \tilde w_k|^2+ (V_\infty - \mu) \tilde w_k^2 +
\phi_{\tilde w_k}\tilde w_k^2
\\
\le \irt |\n \tilde w_k|^2+  V(x)  \tilde w_k^2 + \phi_{\tilde
w_k}\tilde w_k^2 = \irt |\tilde w_k|^{p+1},
\end{multline*}
and then $t_k\le 1$. By \eqref{eq:wk} and ({\bf V3}), we conclude that
\begin{align*}
c(V_\infty -\mu) & \le \frac {t_k^2} 2 \irt |\n \tilde w_k|^2 +
(V_\infty - \mu) \tilde w_k^2 + \frac {t_k^4} 4 \irt \phi_{\tilde
w_k} \tilde w_k^2 - \frac {t_k^{p+1}} {p+1}\irt |\tilde w_k|^{p+1}
\\
& \le \frac {t_k^2} 2 \irt |\n \tilde w_k|^2 + V(x)\tilde w_k^2 +
\frac {t_k^4} 4 \irt \phi_{\tilde w_k} \tilde w_k^2 - \frac
{t_k^{p+1}} {p+1}\irt |\tilde w_k|^{p+1}
\\
& = \left(\frac {t_k^2} 2 - \frac{t_k^{p+1}}{p+1} \right) \irt |\n
\tilde w_k|^2 +  V(x) \tilde w_k^2 + \left(\frac {t_k^4} 4 - \frac
{t_k^{p+1}} {p+1}\right)\irt \phi_{\tilde w_k} \tilde w_k^2
\\
& \le J(\tilde w_k) = J(\eta_k w_k) \le  J(w_k) \le J(u_k) +O(\mu)
= c +o_k(1) +O(\mu),
\end{align*}
but, letting $\mu$ go to zero and $k$ go to $\infty,$ by Lemma
\ref{le:vn}, this yields a contradiction with Lemma
\ref{le:cinfty}.
\\
{\sc Case 3:} up to a subsequence, $G(v_k) > 0$ and $G(w_k) > 0$.
\\
By \eqref{eq:G}, we infer that $G(v_k)=O(\mu)$ and
$G(w_k)=O(\mu)$. Let $(\eta_k)_k$ be such that $\eta_k w_k \in
\Ne$. If $\eta_k\le 1+O(\mu)$, we can repeat the arguments of Case
2. Suppose that
\begin{equation*}
\lim_k \eta_k=\eta_0>1.
\end{equation*}
We have
\begin{align*}
O(\mu)&= G(w_k) = \irt |\n w_k|^2+  V(x) w_k^2 +\phi_{w_k} w_k^2-
|w_k|^{p+1}
\\
&=\left( 1-\frac{1}{\eta_k^{p-1}}\right)\irt |\n w_k|^2+ V(x) w_k^2
+\left( 1-\frac{1}{\eta_k^{p-3}}\right) \irt \phi_{w_k}
w_k^2
\end{align*}
and so
\[
\irt |\n w_k|^2 +  V(x) w_k^2=O(\mu),
\]
which contradicts \eqref{eq:contr2}.
\end{proof}

\begin{proofmain}
By Theorem \ref{th:conc}, for any $\delta>0$ there exists $r>0$
such that
\begin{equation}\label{eq:conc2}
\|u_n\|_{H^1(B_r^c)}<\d,\quad \hbox{uniformly for }n\ge 1.
\end{equation}
By \eqref{eq:weakM} and
\eqref{eq:conc2}, we have that, taken $s\in [2,6[$, for any $\d>0$
there exists $r>0$ such that, for any $n\ge 1$ large enough
\begin{align*}
\|u_n-\bar u\|_{L^s(\RT)} & \le \|u_n-\bar
u\|_{L^s(B_r)}+\|u_n-\bar u\|_{L^s(B^c_r)}
\\
& \le \d + C\left(\|u_n\|_{H^{1}(B_r^c)}+\|\bar u\|_{H^{1}(B_r^c)}\right)\le
(1+2C)\d,
\end{align*}
where $C>0$ is the constant of the embedding
$H^{1}(B_r^c)\hookrightarrow L^s(B^c_r).$ We deduce that
\begin{equation}\label{eq:conv}
u_n\to \bar u\hbox{ in }L^s(\RT),\;\hbox{for any }s\in [2,6[.
\end{equation}
Since $\phi$ is continuous from $L^{12/5}(\RT)$ to $\D$ (see for
instance \cite{Ru}), from \eqref{eq:conv} we deduce that
\begin{align}
\phi_{u_n}  \to \phi_{\bar u} \;\hbox{ in }\D,& \qquad \hbox{as } n\to \infty,  \nonumber\\
\irt\phi_{u_n}u_n^2  \to  \irt\phi_{\bar u}\bar u^2,& \qquad
\hbox{as } n\to \infty,  \label{eq:convphi}
\end{align}
and for any $\psi\in C_0^\infty(\RT)$
\begin{equation}\label{eq:convint}
\irt\phi_{u_n}u_n\psi  \to  \irt\phi_{\bar u}\bar u\psi.
\end{equation}
By \eqref{eq:lim}, we can suppose (see \cite{Wi}) that $(u_n)_n$
is a Palais-Smale sequence for $I|_{\Ne}$ and, as a consequence,
it is easy to see that $(u_n)_n$ is a Palais-Smale sequence for
$I$. By \eqref{eq:weakM}, \eqref{eq:conv} and \eqref{eq:convint},
we conclude that $I'(\bar u)=0.$
\\
Since $(u_n)_n$ is in $\Ne$, by 3 of Lemma \ref{le:N}
$(\|u_n\|_{p+1})_n$ is bounded below by a positive constant. As a
consequence, \eqref{eq:conv} implies that $\bar u\neq 0$ and so
$\bar u\in \Ne$.
\\
Finally, by \eqref{eq:weakM}, \eqref{eq:conv} and
\eqref{eq:convphi} and by ({\bf V2-3}) we get
\begin{align*}
c\le I(\bar u) \le \liminf I(u_n) = c,
\end{align*}
so we can conclude that $(\bar u, \phi_{\bar u})$ is a ground
state solution of \eqref{SM}.
\end{proofmain}

\end{document}